\documentclass{article}
\usepackage{graphicx,amsmath,amssymb,latexsym}
\usepackage{amsfonts}
\usepackage{amsthm}
\usepackage{newlfont}

\newcommand{\Real}{\mathbb R}
\begin{document}

\begin{center}
\textbf{Time and Space Varying Copulas}\vspace{12pt}

Glenis Crane \footnote{Author to whom correspondence should be addressed, Dr G.J. Crane,\\
School of Mathematical Sciences, University of Adelaide, South Australia 5005, Australia, \\
 ph: 618 83036184, Email: glenis.crane@adelaide.edu.au} \\

\textit{School of Mathematical Sciences} \\

\textit{University of Adelaide }\\
\textit{  South Australia}\\
\end{center}
\section{Abstract}
In this article we review existing literature on dynamic copulas and
then propose an $n$-copula which varies in time and space. Our
approach makes use of stochastic differential equations, and
 gives rise to a dynamic copula which is able to capture
the dependence between multiple Markov diffusion processes. This
model is suitable for pricing basket derivatives in finance and may
also be applicable to other areas such as bioinformatics and
environmental science.

\section{Introduction and motivation}
Mapping joint probability distribution functions to copula functions
is straight forward when they are static, due to Sklar's Theorem. On
the other hand, mapping time dependent probability functions to
copula functions is more problematic. In this article we
\begin{description}
\item[(a)] review the techniques for creating time dependent copulas and

\item[(b)] extend the method described in~\cite{GAC06},~\cite{GAH06}, since it
incorporates both time and space. These equations are the first of
their kind in higher dimensions, since only $2$-dimensional examples
have previously been described.
\end{description}

There are at least two areas in which the time dependent copulas of
this chapter are applicable,
\begin{itemize}
\item \textit{Credit derivatives}. We would assume in this application that we have a portfolio of
$n$ firms, and $X_i(t)$ is the value of a $i$-th firm's assets at
time $t$. Each marginal distribution associated with $X_i(t)$ would
represent the probability of the firm's value falling below some
threshold, given certain information at time zero. The time varying
copula would represent the evolution of the joint distribution or
state of the entire portfolio.

\item  \textit{Genetic drift}. For example, each $X_i(t)$ may
represent the frequency of a particular gene at time $t$. Each
marginal distribution would represent the probability that the
frequency of a particular gene had fallen below some threshold. The
copula would relate to the evolution of a group of genes of
interest.
\end{itemize}
\subsection{Notation and Definitions}
In order to understand some of the issues surrounding the mapping of
copulae to distributions it is necessary to go back to some of the
basic definitions and some notation in relation to the probability
distributions of interest.

The notation used for a univariate probability transition function
 here will be
\begin{equation}\label{E:FFX2}
 \textrm{Pr}\{X_i(t_i) \leq x_i \mid X_j(t_j)= x_j\} = F(t_i,x_i
\mid t_j,x_j).
\end{equation}
If $t_j = 0$ then it is quite common to suppress the zero and so the
notation the distribution in this case would be $F(t_i,x_i \mid
x^i_0)$.

\subsection{Method of Darsow et al}
Authors in~\cite{Dar92} were the first to attempt to map a
transition probability function to a copula. First let us recall the
definition
of a bivariate copula.\\

\textbf{Definition 1}. \textit{2-copula}. A function $C:[0,1]^2
\rightarrow [0,1]$ is a copula if it satisfies the following
properties;
\begin{enumerate}
\item  $C(u_1,0) = 0$, $C(0,u_2) = 0$, for all $u_1, u_2 \in
[0,1]$
\item $C(u_1,1) = u_1$, $C(1,u_2) = u_2$, for all $u_1, u_2 \in
[0,1]$ and
\item For every $u_a$, $u_b$, $v_a$, $v_b$ $\in$ [0,1], such
that $u_a \leq u_b$, $v_a \leq v_b$, the volume of $C$,
$V_C([u_a,u_b]\times[v_a,v_b]) \geq 0$, that is
\begin{equation}
\nonumber C(u_b,v_b) - C(u_b,v_a)-C(u_a,v_b) + C(u_a,v_a) \geq 0.
\end{equation}
\end{enumerate}

\textbf{Sklar's Theorem}.  Suppose $H$ is a bivariate joint
distribution with marginal distributions $F_1$ and $F_2$ then there
exists a 2-copula $C$, such that for all $x_1,x_2 \in \bar{\Real}$
\begin{equation}
H(x_1,x_2) = C(F_1(x_1),F_2(x_2)).
\end{equation}
If $F_1$ and $F_2$ are continuous distributions then $C$ is unique,
otherwise $C$ is uniquely determined on $Ran\textsl{$F_1$} \times
Ran\textsl{$F_2$}$, see Nelsen~\cite{Ne99}. \vspace{12pt}

\textbf{Definition 5.1}.  \textit{Markov Property}. A stochastic
process  $X_i(t)$ and $ x_i \in \Real$, $a \leq t \leq b$ is said to
satisfy the Markov property if for any $a \leq t_1 \leq t_2 \ldots,
\leq t_n \leq t$, the equality
\begin{equation}
\nonumber  \textrm{Pr}\{ X_i(t) \leq x_i \mid X_i(t_1),
X_i(t_2),\ldots, X_i(t_n) \} =  \textrm{Pr}\{ X_i(t) \leq x_i \mid
X_i(t_n) \}
\end{equation}
holds for any $x_i \in \Real$. A stochastic process is called a
\textit{Markov Process} if it satisfies the Markov property
described in Definition 5.1.\vspace{12pt}

The following notation will be used for an unconditional probability
function at time  $t_i \geq 0$,
\begin{equation}\label{E:FFX}
 \textrm{Pr}\{X_i(t_i) \leq x_i\} = F_{t_i}(x_i)
\end{equation}
for a stochastic process $X_i(t_i)$  and $x_i \in \Real$. Let
\begin{equation}
\nonumber \nabla_{x_i}F = \frac{\partial F}{\partial x_i},
\end{equation}
then the corresponding density function $f$ in this case is such
that
\begin{equation}
\nonumber f(x_i) =\nabla_{x_i}F(x_i).
\end{equation}
 A univariate transition probability function $F$
(Markov process) can be mapped to a bivariate copula $C$ by setting
\begin{equation}
F(t_i,x_i \mid t_j,x_j) = \nabla_{u_2}
C\left(F_{t_i}(x_i),F_{t_j}(x_j))\right),
\end{equation}
where $\nabla_{u_2}$ is the partial derivative with respect to the
second argument of $C$. This mapping enables us to build in time,
see~\cite{Dar92}. The first marginal distribution is associated with
time $t_i$ and the second with time $t_j$. We take the partial
derivative of the copula, since the probability to which it is
mapped is conditional. This method is particularly useful for
building Markov chains. One of the most important innovations which
enabled the authors in~\cite{Dar92} to link copulas to Markov
processes was to introduce the idea of a copula product;
\vspace{12pt}

\textbf{Definition 5.2}. \textit{Copula product}. Let $C_a$ and
$C_b$ be bivariate copulas, then the product of $C_a$ and $C_b$ is
the function $C_a * C_b : [0,1]^2 \rightarrow [0,1]$, such that
\begin{equation}
(C_a * C_b)(x,y) = \int^1_0 \nabla_z C_a (x,z)\nabla_z C_b (z,y) dz.
\end{equation}

This product is essentially the copula equivalent of the
Chapman-Kolmogorov equation, as stated in Theorem 3.2
of~\cite{Dar92}. We restate that theorem here (with modified notation).\\

\textbf{Theorem 3.2}. Let $X_i(t)$, $t \in T$ be a real stochastic
process, and and for each $s, t \in T$ let $C_{st}$ denote the
copula of the random variables $ X_i(s)$ and $ X_i(t)$. The
following are equivalent: \\
\begin{enumerate}
\item  The transition probabilities $F(t,\mathbf{A} \mid s,x_s) =  \textrm{Pr}\{ X_i(t) \in \mathbf{A}
\mid  X_i(s) = x_s\}$ of the process satisfy the Chapman-Kolmogorov
equations
\begin{equation}
F(t,\mathbf{A} \mid s,x_s) = \int_{\Real} F(t,\mathbf{A} \mid
u,\xi)F(u,d\xi \mid s,x_s)
\end{equation}
for all Borel sets $\mathbf{A}$, for all $s < t \in T$, for all
$u\in (s,t) \cap T$ and for almost all $x_s \in \Real$.
\vspace{12pt}
\item For all $s,u,t \in T$ satisfying $s < u < t$,
\begin{equation}
C_{st} = C_{su}*C_{ut}.
\end{equation}
\end{enumerate}

The work in~\cite{Dar92} has advanced both the theory of copulas and
techniques for building Markov processes. This method is also used
in~\cite{ZIZ05} to formulate a Markov chain model of the dependence
in credit risk. The discrete stochastic variable $X_i(t)$ is
interpreted as the rating grade of a firm at a particular point in
time. A variety of copulas were fitted to the data and gave mixed
results. Therefore, no copula was the best for all data sets. This
type of mapping of the transition distribution to the copula is very
simple, however, one consequence is that an $n$-dimensional
transition function
 requires a $2n$-dimensional copula. In other words, as the dimension of the
copula increases, the calculation of the transition function becomes
more and more computationally cumbersome.\\

The method in~\cite{Dar92} has also been extended in~\cite{SSM03},
so that an $n$-dimensional Markov process can be represented by a
combination of bivariate copulas and margins. Hence,
\begin{eqnarray}
\nonumber &&\hspace*{-15mm} \textrm{Pr}\{X_i(t_1) \leq x_1, \ldots, X_i(t_n) \leq x_n\}  \\
\nonumber &=&\prod^n_{i=2} \textrm{Pr}\{X_i(t_i) \leq x_i \mid
X_i(t_{1}) = x_1,\ldots, X_i(t_{i-1}) = x_{i-1}\}
\textrm{Pr}\{X_i(t_1) \leq x_1\}
\\ \nonumber &=& \prod^n_{i=2} \textrm{Pr}\{X_i(t_i) \leq x_i \mid X_i(t_{i-1}) =
x_{i-1}\} \textrm{Pr}\{X_i(t_1) \leq x_1\}\\
&=&
\frac{\prod^n_{i=2}C_{t_{i-1},t_i}(F_{t_{i-1}}(x_{i-1}),F_{t_{i}}(x_{i}))}{\prod^{n-1}_{i=2}F_{t_{i}}(x_{i})}.
\end{eqnarray}

\subsection{Conditional Copula of Patton}
Another approach to building time into a copula was formulated
in~\cite{PAT01}. In order to explain this approach, we need to
recall more definitions and set up notation. Firstly, let $\mathcal
F$ be a filtration, then
\begin{equation}\label{E:FFC}
 \textrm{Pr}\{X_i \leq x_i \mid \mathcal F\} = F_i(x_i \mid \mathcal F)
\end{equation}
The multivariate analogue of equation \eqref{E:FFC} is
\begin{equation}\label{E:FFC2}
 \textrm{Pr}\{X \leq x \mid \mathcal F\} = H(x \mid \mathcal F),
\end{equation}

for $x = (x_1,x_, \ldots,x_n)^T$ such that the volume of $H$,
$V_H(R) \geq 0$, for all rectangles $R \in \Real^n$ with their
vertices in the domain of $H$, see~\cite{SSM03},

\begin{eqnarray}
\nonumber &&H(+\infty,x_i,+\infty,\ldots, +\infty  \mid \mathcal F)
=F_i(x_i \mid \mathcal F),  \quad \textrm{and}\\
\nonumber &&H(-\infty,x_i,\ldots,x_n \mid \mathcal F) = 0 \quad
\text{for all} \quad x_1,\ldots, x_n \in \Real.
\end{eqnarray}

Here $F_i$ is the $i$-th univariate marginal distribution of $H$.
See~\cite{PAT01} for a bivariate version of $H$. As expected, the
density of the conditional $H$ is
\begin{equation}\label{E:FFC2}
h(x \mid \mathcal F) = \nabla_{x_1,\ldots,x_n}H(x \mid \mathcal F).
\end{equation}

In equation \eqref{E:FFC}, the distribution is atypical since it may
be conditional on a vector of variables, not just one, as opposed to
a typical univariate transition distribution.\\

The author in~\cite{PAT01} mapped the conditional distribution $H(x
\mid \mathcal F)$, defined above, to a copula of the same order.
That is, for all $x_i \in \Real$ and $i= 1,2,\ldots,n$,
\begin{equation}
H(x_1,\ldots,x_n \mid \mathcal F) = C(F_1(x_1 \mid \mathcal
F),F_2(x_2 \mid \mathcal F),\ldots, F_n(x_n \mid \mathcal F)\mid
\mathcal F).
\end{equation}
$\mathcal F$ is a sub-algebra or in other words a conditioning set.
Such conditioning is necessary for $C$ to satisfy all the conditions
of a conventional copula. The relationship between the conditional
density $h$ and copula density $c$ is
\begin{equation}
h(x_1,x_2,\ldots,x_n \mid \mathcal F) = c(u_1,u_2,\ldots, u_n \mid
\mathcal F)\prod^n_{i=1} f_i(x_i \mid \mathcal F),
\end{equation}
where $u_i \equiv F_i(x_i \mid \mathcal F)$, $i=1,2,\ldots, n$ and
$f_i$, $i=1,2,\ldots, n$ are univariate conditional densities.\\

In terms of time varying distributions, we can think of the
conditioning set as the history of all the variables in the
distribution. In the case of the Markov processes, it is only the
last time point which is of importance. The implication of this type
of conditioning is that the marginal distributions in the copula can
no longer be typical transition probabilities, but are atypical
conditional probabilities. Hence, if each $X_i$ represented the
value of an asset at time $t$, the associated distribution $F_i$
would represent the distribution of $X_i$, given that we knew the
value of all the assets in the model, $X_1, X_2,\ldots, X_n$, at
some previous time, for example $t-1$. In other words, we can
rewrite the time-varying version of the distribution and copula
above as
\begin{eqnarray}
\nonumber &&\hspace*{-15mm}H_t(x^1_t,x^2_t,\ldots,x^n_t \mid \mathcal F_{t-1})\\
&=& C_t(F^1_t(x^1_t \mid \mathcal F_{t-1}),F^2_t(x^2_t \mid \mathcal
F_{t-1}),\ldots, F^n_t(x^n_t \mid \mathcal F_{t-1})\mid \mathcal
F_{t-1}),
\end{eqnarray}
where
\begin{equation}
\nonumber \mathcal F_{t-1} =
\sigma(x^1_{t-1},x^2_{t-1},\ldots,x^n_{t-1},x^1_{t-2},x^2_{t-2},\ldots,x^n_{t-2},
\ldots, x^1_1,x^2_1,\ldots,x^n_1).
\end{equation}

In~\cite{PAT01}, the marginal distributions are characterized by
Autoregressive (AR) and generalized autoregressive conditional
heteroskedasticity (GARCH) processes. Ultimately, they are handled
in the same way as other time series processes.

\subsection{Pseudo-copulas of Fermanian and Wegkamp}
As we have seen above, Markov processes are only defined with
respect to their own history, not the history of other processes.
Therefore, the method in~\cite{PAT01} is good for some applications
but not practical for others. If we want marginal distributions of
processes, conditional on their own history, for example Markov
processes, and want to use a mapping similar to that shown
in~\cite{PAT01}, then it is possible via a conditional
pseudo-copula. Authors in~\cite{FEW04} introduced the notion of
conditional pseudo-copula in order to cover a wider range of
applications than the conditional copula in~\cite{PAT01}. The
definition of a pseudo-copula is

\textbf{Definition 5.3}. \textit{Pseudo-copula}. A function
$C:[0,1]^n \rightarrow [0,1]$ is called an $n$-dimensional
pseudo-copula if
\begin{enumerate}
\item  for every $\textbf{u} \in [0,1]^n$, C(\textbf{u}) = 0 when
at least one coordinate of \textbf{u} is zero,

\item $C(1,1,\ldots,1) = 1$, and

\item for every  $\textbf{u}, \textbf{v} \in
[0,1]^n$ such that $\textbf{u} \leq \textbf{v}$, the volume of C,
$V_C \geq 0$.
\end{enumerate}

The pseudo-copula satisfies most of the conditions of a conventional
copula except for $C(1,1,u_k,1,\ldots,1) = u_k$, so the marginal
distributions of a pseudo-copula may not be uniform. The definition
of a conditional pseudo-copula is \\

\textbf{Definition 5.4}. \textit{Conditional pseudo-copula}. Given a
joint distribution $H$ associated with $X_1,X_2,\ldots, X_n$, an
$n$-dimensional conditional
pseudo-copula with respect to sub-algebras \\
$\mathcal F = (\mathcal F_1,\mathcal F_2,\ldots, \mathcal F_n)$ and
$\mathcal G$ is a random function $C(\cdot\mid \mathcal F,\mathcal
G):[0,1]^n \rightarrow [0,1]$ such that
\begin{equation}
H(x_1,x_2,\ldots,x_n) = C(F_1(x_1 \mid \mathcal F_1),F_2(x_2 \mid
\mathcal F_2),\ldots,F_n(x_n \mid \mathcal F_n) \mid \mathcal
F,\mathcal G)
\end{equation}
almost everywhere, for every $(x_1,x_2,\ldots,x_n)^T \in \Real^n$,
 see \cite{FEC04}.

\subsection{Galichon model} More recently, a dynamic bivariate copula was used to correlate Markov
diffusion processes, see~\cite{GAC06},~\cite{GAH06}. Unlike the
previous models of time dependent copulas, this model addresses the
issue of spatial as well as time dependence. The model uses a
partial differential approach to obtain a representation of the time
dependent copula. An outline of the main result follows. Consider
two Markov diffusion processes $X_1(t)$ and $X_2(t)$, $t \in [0,T]$,
which represent two risky financial assets, for example options with
a maturity date $T$. The diffusions are such that
\begin{eqnarray}
\nonumber &&dX_1(t) = \mu_1(X(t))dt + \tilde{\sigma}_1(X(t))dB_1(t) \\
\nonumber &&dX_2(t) = \mu_2(X(t))dt + \tilde{\sigma}_2(X(t))dB_2(t)\\
&&dB_1(t)dB_2(t) = \rho_{12}(X_1(t),X_2(t))dt,
\end{eqnarray}
where $X(t) = (X_1(t),X_2(t))^T$, $\mu_i$, $\tilde{\sigma}_i$, for
$i=1,2$, are the drift and diffusion coefficients, respectively. The
Brownian motion terms are correlated with coefficient $\rho_{12} \in
[-1,1]$. One would like an expression for the evolution of a copula
between the distributions $F_1$, $F_2$ of $X_1(t)$ and $X_2(t)$,
conditional on information at time $t=0$, $\mathcal F_{t_0}$.
Firstly a joint bivariate distribution $H$ is mapped to a copula
$C$, by
\begin{equation}
H(t,x_1,x_2\mid \mathcal F_{t_0}) = C(t,F_1(t,x_1\mid \mathcal
F_{t_0}),F_2(t,x_2\mid \mathcal F_{t_0})\mid\mathcal F_{t_0})
\end{equation}
then the Kolmogorov forward equation is used to obtain an expression
for $\nabla_tC$. Letting
\begin{equation}
\nonumber u_1 = F_1(t,x_1\mid \mathcal F_{t_0})\quad \text{and}
\quad u_2 = F_2(t,x_2\mid \mathcal F_{t_0}), \quad u_1,u_2 \in
[0,1],
\end{equation}
$x = (x_1,x_2)^T$ and shortening the notation for the copula to
$C(t,u_1,u_2)$, then the time dependent copula in~\cite{GAC06} is
\begin{eqnarray}
\nonumber \nabla_tC(t,u_1,u_2)   &=& \frac{1}{2}
\tilde{\sigma}_1^2(x) f_1^2(t,x_1\mid \mathcal
F_{t_0})\nabla^2_{u_1}C(t,u_1,u_2) + \frac{1}{2}
\tilde{\sigma}_2^2(x) f_2^2(t,x_2\mid \mathcal F_{t_0})\nabla^2_{u_2}C(t,u_1,u_2) \\
\nonumber &-& \nabla_{u_1}C(t,u_1,u_2)\mathcal{B}_1F_1(t,x_1\mid
\mathcal F_{t_0}) \\ \nonumber &+&
\int_{(-\infty,x_2]}\hspace*{-2mm}\nabla_{u_1,u_2}C(t,u_1,u_2)f_2(t,z_2\mid
\mathcal F_{t_0})\mathcal{B}_1F_1(t,z_1\mid \mathcal F_{t_0})dz_2\\
\nonumber &-& \nabla_{u_2}C(t,u_1,u_2)\mathcal{B}_2F_2(t,x_2\mid
\mathcal F_{t_0}) \\ \nonumber &+&
\int_{(-\infty,x_1]}\hspace*{-2mm}\nabla_{u_1,u_2}C(t,u_1,u_2)f_1(t,z_1\mid
\mathcal F_{t_0})\mathcal{B}_2F_2(t,z_2\mid \mathcal F_{t_0})dz_1\\
\nonumber &+&
\tilde{\sigma}_1(x)\tilde{\sigma}_2(x)\rho_{12}(x_1,x_2)
f_1(t,x_1\mid \mathcal F_{t_0})f_2(t,x_2\mid \mathcal F_{t_0})
\nabla_{u_1,u_2}C(t,u_1,u_2), \\ &&
\end{eqnarray}
where  $\mathcal{B}_1$ and $\mathcal{B}_2$ are the following
operators, given any function $g \in C^2(\Real)$,
\begin{eqnarray}
\nonumber \mathcal{B}_1 g &=& \left\{
\nabla_{x_1}\bigg(\frac{1}{2}\tilde{\sigma}^2_1(x)\bigg)-\mu_1(x)\right\}\nabla_{x_1}g
+ \bigg(\frac{1}{2}\tilde{\sigma}^2_1(x)\bigg)\nabla^2_{x_1}g \\
\nonumber \mathcal{B}_2 g &=& \left\{
\nabla_{x_2}\bigg(\frac{1}{2}\tilde{\sigma}^2_2(x)\bigg)-\mu_2(x)\right\}\nabla_{x_2}g
+ \bigg(\frac{1}{2}\tilde{\sigma}^2_2(x)\bigg)\nabla^2_{x_2}g
\end{eqnarray}
and
\begin{equation}
\nonumber \nabla_{x_i}g = \frac{\partial g}{\partial x_i}, \quad
\nabla^2_{x_i}g = \frac{\partial^2 g}{\partial x^2_i}.
\end{equation}
For the greatest flexibility we would choose
\begin{equation}
\nonumber inf\{x_i : F_i(t,x_i\mid \mathcal F_{t_0}) \geq u_i\} =
F_i^{-1}(t,u_i \mid \mathcal F_{t_0}), \quad u_i \in [0,1].
\end{equation}
That is, $F_i^{-1}$ is the pseudo-inverse. If $X_1(t)$ and $X_2(t)$
are individually Markov, that is, $\tilde{\sigma}_i$ and $\mu_i$
depend only on $x_i$, for $i = 1,2$, then the formula for the time
dependent copula simplifies to
\begin{eqnarray}
\nonumber && \hspace*{-10mm} \nabla_tC(t,u_1,u_2) =
\tilde{\sigma}_1(x)\tilde{\sigma}_2(x)\rho_{12}(x_1,x_2)
f_1(t,x_1\mid \mathcal F_{t_0})f_2(t,x_2\mid \mathcal
F_{t_0})\nabla_{u_1,u_2}C(t,u_1,u_2)\\  \nonumber &+& \frac{1}{2}
\tilde{\sigma}_1^2(x) f_1^2(t,x_1\mid \mathcal
F_{t_0})\nabla^2_{u_1}C(t,u_1,u_2) + \frac{1}{2}
\tilde{\sigma}_2^2(x) f_2^2(t,x_2\mid \mathcal F_{t_0})\nabla^2_{u_2}C(t,u_1,u_2). \\
&&
\end{eqnarray}

The main aim of this chapter is to extend some of these current
models of time dependent copulas. We derive an $n$-dimensional
version of the model in~\cite{GAC06}. A reformulation is also given,
in which linear combinations of independent Brownian motion terms
are used.

\section{$n$-dimensional Galichon Model for CDOs} Suppose we have an $n \times
n$ system of stochastic differential equations, such that $X(t) \in
\Real^n$  and $B(t)$ is an $n$-dimensional Brownian motion. The
vector $X(t)$ could represent a portfolio of risky assets, as in a
Collateralized Debt Obligation. We want to find a partial
differential equation with respect to a time dependent $n$-copula,
which gives us information on the riskiness of the package of
assets. As in the 2-dimensional model, $t$ is a scalar such that $t
\in (0,T]$. Let $\mathcal{F}_t$ be a $\sigma$-algebra generated by
$\{B(s); s \leq t \}$ and assume $X(t)$ is measurable with respect
to $\mathcal{F}_t$. In this case the diffusions are such that
\begin{eqnarray}\label{E:Gbn}
dX(t) &=& \mu(X(t))dt + \tilde{A}dB(t)
\\ dB_i(t)dB_j(t) &=&
\rho_{ij}(X_i(t),X_j(t))dt,
\end{eqnarray}
where
\begin{equation*}
dX(t) = \left(
\begin{matrix}
dX_1(t) \\
dX_2(t) \\
\vdots \\
dX_n(t)
\end{matrix}
\right), \quad dB(t) = \left(
\begin{matrix}
dB_1(t) \\
dB_2(t) \\
\vdots \\
dB_n(t)
\end{matrix}
\right)
\end{equation*}
\begin{equation*}
\quad \mu(X(t)) = \left(
\begin{matrix}
\mu_1(X(t)) \\
\mu_2(X(t)) \\
\vdots \\
\mu_n(X(t))
\end{matrix}
\right), \quad \tilde{\sigma}(X(t)) = \left(
\begin{matrix}
\tilde{\sigma}_1(X(t)) \\
\tilde{\sigma}_2(X(t)) \\
\vdots \\
\tilde{\sigma}_n(X(t))
\end{matrix}
\right).
\end{equation*}
Note that in this case $\mu$ and $\tilde{\sigma}$ are $n$-vector
functions which represent the drift and diffusion coefficients of
the process, respectively. Let $\tilde{A}$ be
\begin{equation*}
\tilde{A} = diag(\tilde{\sigma}(X(t))) = \left(
\begin{matrix}
\tilde{\sigma}_1(X(t)) & 0 & \ldots & \ldots & 0\\
0 & \tilde{\sigma}_2(X(t)) & 0 & \ldots & 0\\
\vdots & & & & \vdots \\
0 & \ldots & \ldots & 0 & \tilde{\sigma}_n(X(t)). \\
\end{matrix}
\right)
\end{equation*}
The correlation coefficients $\rho_{ij} \in [-1,1]$ and let $\rho$
be
\begin{equation*}
\rho = \left(
\begin{matrix}
1 & \rho_{12}(X_1(t),X_2(t))  & \ldots & \rho_{1n}(X_1(t),X_n(t))\\
\rho_{21}(X_2(t),X_1(t)) & 1  & \ldots & \rho_{2n}(X_2(t),X_n(t))\\
\vdots & &  & \vdots \\
\rho_{n1}(X_n(t),X_1(t)) & \ldots & \ldots & 1 \\
\end{matrix}
\right).
\end{equation*}
\vspace{12pt}

Three conditions are required for the existence and uniqueness of a
solution to equation \eqref{E:Gbn}:
\begin{enumerate}
\item Coefficients $\mu(x)$ and $\tilde{\sigma}(x)$ must be defined for $x \in
\Real^n$ and measurable with respect to $x$.
\item For $x,y \in \Real^n$, there exists a constant $K$
such that \begin{eqnarray}
\nonumber \parallel \mu(x) - \mu(y) \parallel  &\leq& K \parallel x-y \parallel, \\
\nonumber  \parallel \tilde{\sigma}(x) - \tilde{\sigma}(y)
\parallel &\leq& K \parallel x-y \parallel, \\
\nonumber \parallel \mu(x) \parallel^2 + \parallel \tilde{\sigma}(x)
\parallel^2  &\leq& K^2(1 + \parallel x\parallel^2)
\end{eqnarray}
and
\item $X(0)$ does not depend on $B(t)$ and $\mathbb{E}[X(0)^2] < \infty$.
\end{enumerate}
\vspace{12pt}

 \textbf{Theorem 5.1.} The time dependent $n$-copula
$\nabla_tC(t,u)$ between a vector of distributions $u_i =
F_i(t,x_i\mid x_0)$, $i=1,\ldots,n$, associated with the Markov
diffusions $X(t) = [X_1(t),\ldots,X_n(t)]^T$,  conditional on
information at time $t=0$, $\mathcal{F}_{t_0} = x_0$ is
\begin{eqnarray}
\nonumber &&\hspace*{-8mm}\nabla_tC(t,u) = \frac{1}{2}
\sum^n_{i=1}\int_{(-\infty,\bar{x}]}\hspace*{-2mm}
\tilde{\sigma}_i(z)^2f^2_i(t,x_i\mid
x_0)\nabla_{z_1,.,\hat{z}_i,.,z_n}\nabla_{u_i}^2C(t,u)
d\bar{z}\\
\nonumber \hspace*{-1mm}&+& \hspace*{-2mm}
\sum^n_{i=1}\bigg(\hspace*{-2mm}- \nabla_{u_i}C(t,u)\mathcal{B}^i_t
F_i(t,x_i\mid x_0) +\int_{(-\infty,\bar{x}]}
\hspace*{-5mm}\nabla_{z_1,.,\hat{z}_i,.,z_n}\hspace*{-2mm}\nabla_{u_i}C(t,u)\mathcal{B}^i_t
F_i(t,z_i\mid x_0) d\bar{z}\bigg) \\
\nonumber \hspace*{-5mm}&+&
\hspace*{-3mm}\frac{1}{2}\sum^n_{\substack{i,j=1\\ i \neq j}}
\hspace*{-1mm}
\int_{(-\infty,\check{x}]}\hspace*{-9mm}\rho_{ij}(x_i,x_j)\tilde{\sigma}_i(z)\tilde{\sigma}_j(z)f_i
(t,x_i \hspace*{-1mm}\mid \hspace*{-1mm} x_0)f_j(t,x_j
\hspace*{-1mm}\mid \hspace*{-1mm}
x_0)\nabla_{z_1,.,\hat{z}_i,\hat{z}_j,.,z_n}\hspace*{-1mm}\nabla_{u_i,u_j}C(t,u)
d\check{z}, \\
&&
\end{eqnarray}
where $x_i = F_i^{-1}(t,u_i \mid \mathcal{F}_{t_0})$,
$i=1,\ldots,n$. The intervals for the integration are
\begin{eqnarray}
\nonumber (-\infty,\bar{x}] &=&
(-\infty,x_1]\times\ldots\times(-\infty,x_{i-1}]\times
(-\infty,x_{i+1}]\times \ldots \times(-\infty,x_n] \quad \text{and} \quad\\
\nonumber (-\infty,\check{x}] &=&
(-\infty,x_1]\times\ldots\times(-\infty,x_{i-1}]\times
(-\infty,x_{i+1}]\times \ldots\times(-\infty,x_{j-1}] \\
\nonumber &&\times(-\infty,x_{j+1}]\ldots\times(-\infty,x_n].
\end{eqnarray}
 Also note that
\begin{eqnarray}
\nonumber d\bar{z} &=& dz_1dz_2\ldots
dz_{i-1}dz_{i+1}\ldots dz_{n-1}dz_n \quad \text{and} \quad  \\
\nonumber d\check{z} &=& dz_1dz_2 \ldots dz_{i-1}dz_{i+1}\ldots
dz_{j-1}dz_{j+1}\ldots dz_{n-1}dz_n.
\end{eqnarray}
Thus, the $i$-th term is excluded in the first two integrals on the
right hand side of Theorem 5.1. Similarly, in the last integral the
$i$-th and $j$-th terms are excluded. Furthermore, for any smooth
function $g$
\begin{eqnarray}
\nonumber \nabla_{z_1,.,\hat{z}_i,.,z_n} g &=& \frac{\partial
^{n-1}g}{\partial z_1 \ldots\partial z_{i-1}\partial
z_{i+1}\ldots\partial z_{n}} \\
\nonumber \text{and} \quad \nonumber
\nabla_{z_1,.,\hat{z}_i,\hat{z}_j,.,z_n} g &=& \frac{\partial
^{n-2}g}{\partial z_1 \ldots\partial z_{i-1}\partial z_{i+1}\ldots
z_{j-1}\partial z_{j+1}\ldots\partial z_{n}}.
\end{eqnarray}
 The operators
$\mathcal{B}_t^i$, $i=1,\ldots, n$ are the same as in the two
dimensional model. If the diffusions are individually Markov, that
is, each $\sigma_k$ and $\mu_k$ depends only on $x_k$ then
 the expression for $\nabla_tC$ simplifies to
\begin{eqnarray}
\nonumber &&\hspace*{-10mm}\nabla_tC(t,u) = \frac{1}{2}\sum^n_{i=1}
\tilde{\sigma}_i(x_i)^2 f^2_i(t,x_i\mid x_0) \nabla^2_{u_i}C(t,u)
\\ \nonumber &+& \hspace*{-2mm}\frac{1}{2}\text{Tr} \big\{\big[\mathcal{H}_u^C(t,u)-
\text{diag}\{\nabla_{u_1}^2 C(t,u),\nabla_{u_2}^2
C(t,u),\ldots,\nabla_{u_n}^2
C(t,u)\}\big]D\tilde{A}\rho\tilde{A}^TD^T\big\},
\end{eqnarray}
where
\begin{equation*}
D = \left(
\begin{matrix}
f_1 & 0 & \hdots & 0\\
0 & f_2 & \hdots & 0\\
0  & 0  &  \hdots   & 0 \\
0 & 0 & \hdots & f_n
\end{matrix}
\right).
\end{equation*}

\vspace{12pt} \textbf{Proof.}
 In this case 1-dimensional Ito formula for each
component of $X(t)$ is
\begin{eqnarray}
\nonumber dg(X_i(t)) &=& \big\{\nabla_{x_i}g(X_i(t))\mu_i(X(t)) +
\frac{1}{2} \nabla^2_{x_i}g(X_i(t)) \tilde{\sigma}_i^2(X(t))\big\}dt
\\ \nonumber &+& \nabla_{x_i}g(X_i(t))\tilde{\sigma}_i(X(t))
dB_i(t).
\end{eqnarray}
Define the vector $\nabla_x$ of partial derivatives with respect to
components of $x$, as
\begin{equation*}
\nabla_x g(X(t)) = \left(
\begin{matrix}
\nabla_{x_1}g(X(t)) \\
\nabla_{x_2}g(X(t)) \\
\vdots \\
\nabla_{x_n}g(X(t))
\end{matrix}
\right)
\end{equation*}
and the Hessian matrix of $g(X(t))$
\begin{equation}
\mathcal{H}_x^g(X(t)) \equiv \bigg( \bigg(  \nabla_{x_i
x_j}g(X(t))\bigg)\bigg)_{1 \leq i,j \leq n }.
\end{equation}
In this case, assume $g \in C^2(\Real^n)$, then the $n$-dimensional
Ito formula for $g(X(t))$ is
\begin{eqnarray}
\nonumber  &&\hspace*{-8mm}dg(X(t)) \\
&=& \big\{\langle \nabla_x g(X(t)),\mu(X(t)) \rangle + \frac{1}{2}
\text{Tr} \big(\mathcal{H}^g_x(X(t))\tilde{A}\rho\tilde{A}^T
\big)\big\}dt \nonumber + \nabla_xg(X(t))^T\tilde{A}dB(t),
\end{eqnarray}

where $\langle a,b \rangle =a^Tb$ for any vectors $a$ and $b$. Let
the operators $\mathcal{A}$ on distributions (Kolmogorov backward
equations), analogous to those in~\cite{GAC06},~\cite{GAH06}, be
called $\mathcal{A}^i_t$ and $\mathcal{A}^n_t$ for the 1- an
$n$-dimensional case, respectively. With respect to typical
distributions $F_i(t,x_i\mid \tau,\xi_i)$ and $H(t,x\mid \tau,\xi)$,
the operators are
\begin{equation}
\mathcal{A}^i_tF_i(t,x_i\mid \tau,\xi_i) =
\mu_i(x)\nabla_{\xi_i}F_i(t,x_i\mid \tau,\xi_i) +
\frac{1}{2}\tilde{\sigma}_i^2 \nabla_{\xi_i}^2F_i(t,x_i\mid
\tau,\xi_i)
\end{equation}
and
\begin{equation}
 \mathcal{A}^n_tH(t,x\mid \tau,\xi) = \langle \nabla_{\xi}
H(t,x\mid \tau,\xi),\mu(x) \rangle + \frac{1}{2} \text{Tr}
\left(\mathcal{H}^H_{\xi}(t,x\mid \tau,\xi)\tilde{A}\rho\tilde{A}^T
\right).
\end{equation}

The operators $\mathcal{A}^i_t$, $i=1,\ldots,n$ and
$\mathcal{A}^n_t$ are not used in the rest of the formulation, but
are mentioned briefly, in view of the fact the Kolmogorov forward
equations, which are required, are the associated adjoint operators
of these. Assuming the density functions of $H$ and $F$ are $h$ and
$f$, respectively, then the adjoint operators $\mathcal{A}^{i*}_t$,
$i=1,\ldots,n$ and $\mathcal{A}^{n*}_t$ have the form
\begin{equation}
\mathcal{A}^{i*}_tf_i(t,x_i\mid \tau,\xi_i) = -
\nabla_{x_i}\big[\mu_i(x)f_i(t,x_i\mid \tau,\xi_i)\big] +
\nabla_{x_i}^2\big[\frac{1}{2}\tilde{\sigma}_i^2f_i(t,x_i\mid
\tau,\xi_i)\big]
\end{equation}
and
\begin{eqnarray}
\nonumber &&\hspace*{-10mm}\mathcal{A}^{n*}_th(t,x\mid \tau,\xi)\\
\nonumber &=& \hspace*{-2mm}- \sum_{i=1}^n
\hspace*{-1mm}\nabla_{x_i}\big[\mu_i(x)h(t,x\mid \tau,\xi) \big]
+\frac{1}{2}\sum_{i,j=1}^n
\hspace*{-1mm}\nabla_{x_i,x_j}\big[\rho_{ij}(x_i,x_j)\tilde{\sigma}_i(x)\tilde{\sigma}_j(x)h(t,x\mid
\tau,\xi) \big].\\
&&
\end{eqnarray}
The marginal density functions $f_i$ and joint density $h$, are such
that \\
$f_i(t,x_i\mid \mathcal{F}_{t_0}) = f_i(t,x_i\mid x_0)$, and
$H(t,x\mid \mathcal{F}_{t_0}) = H(t,x\mid x_0)$, where \\
$x_0 = (x_1 = X_1(0),x_2 = X_2(0),\ldots,x_n = X_n(0))$, see
Appendix 5.A. In other words, the assumption made here is that all
the distributions are conditional on the entire vector of
realizations of $x$ at time zero. As in the 2-dimensional case, it
is possible to express the operator $\mathcal{A}^{n*}_t$ in terms of
the operators $\mathcal{A}^{i*}_t$ associated with the univariate
distributions;
\begin{equation}
\mathcal{A}^{n*}_tg= \sum_{i=1}^n \mathcal{A}^{i*}_tg +
\frac{1}{2}\sum^n_{\substack{i,j=1 \\ i \neq j}}\nabla_{x_i,x_j}
\big[\rho_{ij}(x_i,x_j)\tilde{\sigma}_i(x)\tilde{\sigma}_j(x)
g\big].
\end{equation}
Given that
\begin{equation}\label{E:lrf}
\nabla_tf_i(t,x_i\mid x_0) = \mathcal{A}^{i*}_tf_i(t,x_i\mid x_0),
\end{equation}
we can integrate the left hand side of \eqref{E:lrf} with respect to
$x_i$, call it $\mathcal{B}^i_t$, and we obtain
\begin{eqnarray}
\nonumber \mathcal{B}^i_tF_i(t,x_i\mid x_0) &=&
\int_{(-\infty,x_i]}\hspace*{-2mm}\nabla_tf_i(t,z_i\mid x_0)dz_i \\
\nonumber &=& \int_{(-\infty,x_i]}\hspace*{-2mm}\nabla_t\nabla_{z_i}F_i(t,z_i\mid x_0)dz_i \\
&=& \nabla_tF_i(t,x_i\mid x_0).
\end{eqnarray}
Integrating the right hand side of \eqref{E:lrf} with respect to
$x_i$ gives us
\begin{multline}
\nonumber
\int_{(-\infty,x_i]}\hspace*{-2mm}\mathcal{A}^{i*}_tf_i(t,z_i\mid
x_0)dz_i  = -\mu_i(x)f_i(t,x_i\mid x_0) +
\nabla_{x_i}\big[\frac{1}{2}\tilde{\sigma}^2_i(x)f_i(t,x_i\mid
x_0)\big]\\
=\big[\nabla_{x_i}\{\frac{1}{2}\tilde{\sigma}^2_i(x)\}-\mu_i(x)\big]\nabla_{x_i}F_i(t,x_i\mid
x_0) + \frac{1}{2}\tilde{\sigma}^2_i(x)\nabla^2_{x_i}F_i(t,x_i\mid
x_0),
\end{multline}
so
\begin{equation}
\mathcal{B}^i_tF_i(t,x_i\mid x_0) =
\big[\nabla_{x_i}\{\frac{1}{2}\tilde{\sigma}^2_i(x)\}-\mu_i(x)\big]\nabla_{x_i}F_i(t,x_i\mid
x_0) + \frac{1}{2}\tilde{\sigma}^2_i(x)\nabla^2_{x_i}F_i(t,x_i\mid
x_0).
\end{equation}
Similarly, integrating over $\mathcal{A}^{n*}_t$ will give us the
analogous operator $\mathcal{B}^i_n$ for the multivariate
distribution $H$. Now, since
\begin{eqnarray}\label{E:bb1}
\nonumber&& \hspace*{-7mm}
\int_{(-\infty,x]}\hspace*{-1mm}\mathcal{A}^{n*}_th(t,z\mid x_0)dz
\nonumber \\
&&\hspace*{-5mm}=\sum_{i=1}^n
\int_{(-\infty,x]}\hspace*{-6mm}\mathcal{A}^{i*}_th(t,z\mid x_0)dz
\nonumber + \frac{1}{2}\sum^n_{\substack{i,j=1 \\ i \neq
j}}\int_{(-\infty,x]}\hspace*{-6mm}\nabla_{z_i,z_j}
\big[\rho_{ij}(z_i,z_j)\tilde{\sigma}_i(z)\tilde{\sigma}_j(z)
h(t,z\mid x_0)\big]dz, \\
&&
\end{eqnarray}
where $(-\infty,x] = (-\infty,x_1]\times \ldots\times(-\infty,x_n]$,
it is possible to get an expression for $\mathcal{B}^n_t$ in terms
of $\mathcal{B}^i_t$. That is, let $\mathcal{B}^n_t H(t,x\mid x_0) =
\nabla_t H(t,x\mid x_0)$ and given that $h(t,x\mid x_0) =
\nabla_{x_1,.,x_n}H(t,x\mid x_0)$, we have
\begin{eqnarray}\label{E:bb2}
\nonumber \mathcal{B}^n_t H(t,x\mid x_0) &=&
\frac{1}{2}\sum^n_{\substack{i,j=1 \\ i \neq
j}}\int_{(-\infty,x]}\hspace*{-2mm}\nabla_{z_i,z_j}
\big[\rho_{ij}(z_i,z_j)\tilde{\sigma}_i(z)\tilde{\sigma}_j(z)
\nabla_{z_1,\ldots,z_n}H(t,z\mid x_0)\big]dz\\
&+& \sum_{i=1}^n
\int_{(-\infty,x]}\hspace*{-1mm}\mathcal{A}^{i*}_t\nabla_{z_1,\ldots,z_n}H(t,z\mid
x_0)dz.
\end{eqnarray}
The right hand side of equation \eqref{E:bb2} can be expressed in
terms in terms of the univariate operators $\mathcal{B}^i_t$,
$i=1,2\ldots,n$.
\begin{eqnarray}\label{E:bb3}
\nonumber \mathcal{B}^n_t H(t,x\mid x_0) &=& \hspace*{-2mm}
\frac{1}{2}\sum^n_{\substack{i,j=1 \\ i \neq
j}}\int_{(-\infty,\check{x}]}\hspace*{-4mm}\rho_{ij}(x_i,x_j)\tilde{\sigma}_i(z)\tilde{\sigma}_j(z)\nabla_{z_1,.,\hat{z}_i,\hat{z}_j,.,z_n}\nabla_{z_i,z_j}H(t,z\mid
x_0) d\check{z}\\
 &+& \sum_{i=1}^n
\int_{(-\infty,\bar{x}]}\hspace*{-1mm}\mathcal{B}^i_t\nabla_{z_1,.,\hat{z}_i,.,z_n}H(t,z\mid
x_0)d\bar{z}.
\end{eqnarray}
Let
\begin{equation}
H(t,x\mid x_0) =C(t,F_1(t,x_1\mid x_0),F_2(t,x_2\mid
x_0),\ldots,F_n(t,x_n\mid x_0)\mid x_0)
\end{equation}
where $C$ is an $n$-copula defined on $[0,T]\times[0,1]^n$. At this
point we shorten the notation so that $C(t,F(t,x\mid x_0))$ is the
same copula as above. We now seek an expression for
$\mathcal{B}^n_tC(t,F(t,x\mid x_0))$ by substituting for $H$ with
$C$ in equation \eqref{E:bb3}. Letting $F_i(t,x_i\mid x_0) = u_i$,
 $i =1,2,\ldots,n$, and $u=(u_1,\ldots,u_n)^T$, then from the
first term in equation \eqref{E:bb3} we obtain (see overpage)

\begin{eqnarray}
\nonumber
&&\hspace*{-15mm}\sum_{i=1}^n\int_{(-\infty,\bar{x}]}\hspace*{-1mm}\mathcal{B}^i_t\nabla_{z_1,.,\hat{z}_i,.,z_n}H(t,z\mid
x_0)d\bar{z} \\ \nonumber &=&
\sum_{i=1}^n\int_{(-\infty,\bar{x}]}\hspace*{-1mm}\mathcal{B}^i_t\nabla_{z_1,.,\hat{z}_i,.,z_n}C(t,F(t,z\mid
x_0))d\bar{z} \\ \nonumber &=& \sum_{i=1}^n
\bigg(\int_{(-\infty,\bar{x}]}\big\{\nabla_{z_i}\frac{\tilde{\sigma}_i^2(z)}{2}
- \mu_i(z) \big\} \nabla_{z_i}\nabla_{z_1,.,\hat{z}_i,.,z_n}
C(t,F(t,z\mid x_0))d\bar{z} \\ \nonumber &+&
\int_{(-\infty,\bar{x}]}\frac{\tilde{\sigma}_i^2(z)}{2}
\nabla_{z_i}^2\nabla_{z_1,.,\hat{z}_i,.,z_n} C(t,F(t,z\mid
x_0))d\bar{z} \bigg)\\
 \nonumber &=& \sum_{i=1}^n \bigg(\int_{(-\infty,\bar{x}]}\big\{\nabla_{z_i}\frac{\tilde{\sigma}_i^2(z)}{2}
- \mu_i(z) \big\} f_i(t,z_i\mid
x_0)\nabla_{u_i}\nabla_{z_1,.,\hat{z}_i,.,z_n} C(t,u)d\bar{z} \\
\nonumber &+&
\int_{(-\infty,\bar{x}]}\frac{\tilde{\sigma}_i^2(z)}{2}f_i^2(t,z_i\mid
x_0) \nabla_{u_i}^2\nabla_{z_1,.,\hat{z}_i,.,z_n} C(t,u)d\bar{z} \\
\nonumber &+&
\int_{(-\infty,\bar{x}]}\frac{\tilde{\sigma}_i^2(z)}{2}\nabla_{z_i}f_i(t,z_i\mid
x_0) \nabla_{u_i}\nabla_{z_1,.,\hat{z}_i,.,z_n} C(t,u)d\bar{z} \bigg) \\
 \nonumber &=&  \sum_{i=1}^n \bigg(\int_{(-\infty,\bar{x}]}\big\{\nabla_{z_i}\frac{\tilde{\sigma}_i^2(z)}{2}
- \mu_i(z) \big\} \nabla_{z_i}F_i(t,z_i\mid
x_0)\nabla_{u_i}\nabla_{z_1,.,\hat{z}_i,.,z_n} C(t,u)d\bar{z} \\
\nonumber &+&
\int_{(-\infty,\bar{x}]}\frac{\tilde{\sigma}_i^2(z)}{2}\nabla_{z_i}^2F_i(t,z_i\mid
x_0) \nabla_{u_i}\nabla_{z_1,.,\hat{z}_i,.,z_n} C(t,u)d\bar{z} \\
\nonumber &+&
\int_{(-\infty,\bar{x}]}\frac{\tilde{\sigma}_i^2(z)}{2}f_i^2(t,z_i\mid
x_0) \nabla_{u_i}^2\nabla_{z_1,.,\hat{z}_i,.,z_n} C(t,u)d\bar{z} \bigg) \\
\nonumber &=&
\sum_{i=1}^n\int_{(-\infty,\bar{x}]}\nabla_{z_1,.,\hat{z}_i,.,z_n}
\nabla_{u_i}C(t,u)\mathcal{B}^i_tF_i(t,z_i\mid x_0)d\bar{z}\\
\nonumber &+& \frac{1}{2}\sum_{i=1}^n
\int_{(-\infty,\bar{x}]}\tilde{\sigma}_i^2(z)f_i^2(t,z_i\mid x_0)
\nabla_{z_1,.,\hat{z}_i,.,z_n}\nabla_{u_i}^2 C(t,u)d\bar{z}.
\end{eqnarray}
Since $z$ is a dummy variable  and the multiple integrals exclude
that over $(-\infty, x_i]$, we can write
\begin{eqnarray}
\nonumber
&&\hspace*{-15mm}\sum_{i=1}^n\int_{(-\infty,\bar{x}]}\hspace*{-1mm}\mathcal{B}^i_t\nabla_{z_1,.,\hat{z}_i,.,z_n}H(t,z\mid
x_0)d\bar{z} \\
\nonumber &=&
\sum_{i=1}^n\int_{(-\infty,\bar{x}]}\nabla_{z_1,.,\hat{z}_i,.,z_n}
\nabla_{u_i}C(t,u)\mathcal{B}^i_tF_i(t,z_i\mid x_0)d\bar{z}\\
 &+& \frac{1}{2}\sum_{i=1}^n
\int_{(-\infty,\bar{x}]}\tilde{\sigma}_i^2(z)f_i^2(t,x_i\mid x_0)
\nabla_{z_1,.,\hat{z}_i,.,z_n}\nabla_{u_i}^2 C(t,u)d\bar{z}.
\end{eqnarray}
 From the second
term in equation \eqref{E:bb3} we have
\begin{eqnarray}
\nonumber &&\hspace*{-10mm}\frac{1}{2}\sum^n_{\substack{i,j=1 \\ i
\neq
j}}\int_{(-\infty,\check{x}]}\hspace*{-2mm}\rho_{ij}(x_i,x_j)\tilde{\sigma}_i(z)\tilde{\sigma}_j(z)\nabla_{z_1,.,\hat{z}_i,\hat{z}_j,.,z_n}\hspace*{-1mm}\nabla_{z_i,z_j}H(t,z\mid
x_0) d\check{z} \\
\nonumber &=& \frac{1}{2}\sum^n_{\substack{i,j=1 \\ i \neq
j}}\int_{(-\infty,\check{x}]}\hspace*{-2mm}\rho_{ij}(x_i,x_j)\tilde{\sigma}_i(z)\tilde{\sigma}_j(z)\nabla_{z_1,.,\hat{z}_i,\hat{z}_j,.,z_n}\hspace*{-1mm}\nabla_{z_i,z_j}C(t,F(t,z\mid
x_0)) d\check{z} \\
\nonumber \hspace*{-3mm}&=&
\hspace*{-2mm}\frac{1}{2}\sum^n_{\substack{i,j=1 \\ i \neq
j}}\int_{(-\infty,\check{x}]}\hspace*{-9mm}\rho_{ij}(x_i,x_j)\tilde{\sigma}_i(z)\tilde{\sigma}_j(z)f_i(t,x_i\mid
\hspace*{-1mm} x_0)f_j(t,x_j\mid \hspace*{-1mm}
x_0)\nabla_{z_1,.,\hat{z}_i,\hat{z}_j,.,z_n}\hspace*{-1mm}\nabla_{u_i,u_j}\hspace*{-1mm}C(t,u)
d\check{z}\\
&&
\end{eqnarray}
so
\begin{eqnarray}\label{E:Bn2}
\nonumber &&\hspace*{-10mm}\mathcal{B}^n_tC(t,u)   =
\sum_{i=1}^n\int_{(-\infty,\bar{x}]}\nabla_{z_1,.,\hat{z}_i,.,z_n}
\nabla_{u_i}C(t,u)\mathcal{B}^i_tF_i(t,z_i\mid x_0)d\bar{z}\\
\nonumber &+& \frac{1}{2}\sum_{i=1}^n
\int_{(-\infty,\bar{x}]}\tilde{\sigma}_i^2(z)f_i^2(t,x_i\mid x_0)
\nabla_{z_1,.,\hat{z}_i,.,z_n}\nabla_{u_i}^2 C(t,u)d\bar{z}\\
\hspace*{-3mm}&+& \hspace*{-2mm}\nonumber
\frac{1}{2}\sum^n_{\substack{i,j=1 \\ i \neq
j}}\int_{(-\infty,\check{x}]}\hspace*{-9mm}\rho_{ij}(x_i,x_j)\tilde{\sigma}_i(z)\tilde{\sigma}_j(z)f_i(t,x_i\mid
\hspace*{-1mm} x_0)f_j(t,x_j\mid \hspace*{-1mm}
x_0)\nabla_{z_1,.,\hat{z}_i,\hat{z}_j,.,z_n}\hspace*{-1mm}\nabla_{u_i,u_j}\hspace*{-1mm}C(t,u)
d\check{z}. \\
&&
\end{eqnarray}
Now, we also have
\begin{eqnarray}\label{E:Bn3}
\nonumber \nabla_tH(t,x\mid x_0)  &=& \mathcal{B}^n_tH(t,x\mid x_0) \\
\nonumber &=& \nabla_tC(t,F(t,x\mid x_0)) + \sum^n_{i=1}
\nabla_{u_i}C(t,F(t,x\mid x_0))  \nabla_t F_i(t,x_i\mid x_0) \\
&=& \nabla_tC(t,u) + \sum^n_{i=1} \nabla_{u_i}C(t,u)\mathcal{B}^i_t
F_i(t,x_i\mid x_0).
\end{eqnarray}
Matching equation $\eqref{E:Bn2}$ and $\eqref{E:Bn3}$ and
rearranging, we obtain
\begin{eqnarray}
\nonumber && \hspace*{-10mm}\nabla_tC(t,u) \\
\nonumber &=& \hspace*{-2mm} \sum^n_{i=1}\bigg(\hspace*{-2mm}-
\nabla_{u_i}C(t,u)\mathcal{B}^i_t F_i(t,x_i\mid x_0) +
\int_{(-\infty,\bar{x}]}\hspace*{-7mm}\nabla_{z_1,.,\hat{z}_i,.,z_n}
\nabla_{u_i}C(t,u)\mathcal{B}^i_tF_i(t,z_i\mid x_0)d\bar{z} \bigg) \\
\nonumber &+& \frac{1}{2}\sum_{i=1}^n
\int_{(-\infty,\bar{x}]}\hspace*{-3mm}\tilde{\sigma}_i^2(z)f_i^2(t,x_i\mid
x_0) \nabla_{z_1,.,\hat{z}_i,.,z_n}\nabla_{u_i}^2 C(t,u)d\bar{z}\\
&+& \nonumber \hspace*{-3mm}\frac{1}{2}\sum^n_{\substack{i,j=1 \\ i
\neq j}}
\hspace*{-1mm}\int_{(-\infty,\check{x}]}\hspace*{-9mm}\rho_{ij}(x_i,x_j)\tilde{\sigma}_i(z)\tilde{\sigma}_j(z)f_i(t,x_i
\hspace*{-1mm}\mid \hspace*{-1mm} x_0)f_j(t,x_j \hspace*{-1mm}\mid
\hspace*{-1mm}
x_0)\nabla_{z_1,.,\hat{z}_i,\hat{z}_j,.,z_n}\hspace*{-1mm}\nabla_{u_i,u_j}C(t,u)
d\check{z}.
\end{eqnarray}
If the equations are individually Markov, so that each $\sigma_k$
and $\mu_k$ depends only on $x_k$, then
\begin{equation}
\nonumber \sum^n_{i=1}\bigg(\hspace*{-2mm}-
\nabla_{u_i}C(t,u)\mathcal{B}^i_t F_i(t,x_i\mid x_0) +
\int_{(-\infty,\bar{x}]}\hspace*{-5mm}\nabla_{z_1,.,\hat{z}_i,.,z_n}
\nabla_{u_i}C(t,u)\mathcal{B}^i_tF_i(t,z_i\mid x_0)d\bar{z} \bigg) =
0,
\end{equation}
so the expression for $\nabla_tC$ simplifies to
\begin{eqnarray}
\nonumber && \hspace*{-10mm}\nabla_tC(t,u) = \frac{1}{2}\sum^n_{i=1}
\tilde{\sigma}_i(x_i)^2 f^2_i(t,x_i\mid x_0) \nabla^2_{u_i}C(t,u)
\\ \nonumber &+& \hspace*{-3mm}\frac{1}{2}\text{Tr} \big\{\big[\mathcal{H}_u^C(t,u)-
\text{diag}\{\nabla_{u_1}^2 C(t,u),\nabla_{u_2}^2
C(t,u),\ldots,\nabla_{u_n}^2
C(t,u)\}\big]D\tilde{A}\rho\tilde{A}^TD^T\big\}.
\end{eqnarray}
{\hfill $\square$}
\section{Conclusion}
We have described a dynamic $n$-copula which varies in time and
space. This copula is the first of its kind in greater than 2
dimensions. The dynamic 2-copula was previously described in
\cite{GAH06}. In that case, the copula could be applied the pricing
of pairs of options and other credit derivatives. In the
$n$-dimensional case, it is possible to use the dynamic copula for
the pricing of any basket derivatives or a number of commodities.
Future work in this area may involve numerical experiments,
sensitivity testing and simulations in order to determine how robust
the model is. Other possible applications include that of the health
industry and environmental science.\vspace{12pt}

I would like to acknowledge Alfred Galichon for his valuable
discussions in relation to this work.
\bibliographystyle{amsplain}

\bibliography{dyno}

\end{document}